\documentclass[11pt]{amsart}

\usepackage{url}
\usepackage{graphicx}
\usepackage{color}
\usepackage{tikz}
\usepackage{tikz-3dplot}
\usepackage[labelformat=empty]{subfig}
\usepackage[euler-digits]{eulervm}
\usepackage[all]{xy}

\author[D.Karp]{Dagan Karp}
\email{dagan.karp@hmc.edu}
\author[D.Ranganathan]{Dhruv Ranganathan}
\email{dhruv\_ranganathan@hmc.edu}
\author[P.Riggins]{Paul Riggins}
\email{paul\_riggins@hmc.edu}
\author[U.Whitcher]{Ursula Whitcher}
\email{whitchua@uwec.edu}

\address{Department of Mathematics, Harvey Mudd College}
\address{Department of Mathematics, University of Wisconsin-Eau Claire}

\newtheorem{theorem}{Theorem}
\newtheorem{corollary}[theorem]{Corollary}

\newtheorem*{rem1}{Remark}
\newenvironment{remark}{\begin{rem1}\em}{\end{rem1}}

\newcommand{\CC} {{\mathbb C}}

\newcommand{\PP}{\mathbb{P}}         
\newcommand{\QQ} {{\mathbb Q}}		
		
\newcommand{\ZZ} {{\mathbb Z}}		

\DeclareMathOperator{\ev}{ev}

\DeclareMathOperator{\Aut}{Aut}

\DeclareMathOperator{\vir}{vir}

\newcommand{\cal}{\mathcal}

\def\cI{{\cal I}}

\def\cM{{\cal M}}

\def\cO{{\cal O}}

\def\fI{\mathfrak{I}}

\newcommand{\Mbar}{\overline{\cM}}

\begin{document}

\title{Toric Symmetry of $\CC  \PP^{3}$ }

\begin{abstract} We exhaustively analyze the toric symmetries of $\CC \PP^{3}$ and its toric blowups. Our motivation is to study {\emph{toric symmetry}} as a computational technique in Gromov-Witten theory and Donaldson-Thomas theory. We identify all nontrivial toric symmetries. The induced nontrivial isomorphisms lift and provide new symmetries at the level of Gromov-Witten Theory and Donaldson-Thomas Theory. The polytopes of the toric varieties in question include the permutohedron, the cyclohedron, the associahedron, and in fact all graph associahedra, among others.
\end{abstract}

\maketitle

\section{Introduction}\label{sec: intro}

The group $\Aut (X)$ of automorphisms of a toric variety $X$ is generated by three classes of automorphism: the automorphisms extending from automorphisms of the torus itself, automorphisms corresponding to the roots of $X$, and toric symmetries of $X$. A {\emph{toric symmetry}} of $X$ is simply an automorphism of $X$ induced by an automorphism of the fan $\Sigma_{X}$ (or polytope $\triangle_X$) of $X$. There are many sources for this beautiful material; see, for example~\cite[Theorem 3.6.1]{CoxKatz}.  

How can we analyze the toric symmetries of $X$? If $\tau :X \rightarrow X$ is a toric symmetry of $X$, then $\tau$ may be the identity on homology, or $\tau$ may correspond to a simple relabeling of the cones in $\Sigma_{X}$. Which toric symmetries act on Chow in a manner meaningfully different from the identity?  We ask this more delicate question, and search for toric symmetries with interesting action on homology. We describe this condition in greater detail below, and a complete discussion is found in Section~\ref{sec: toric_symmetries}.

In this work we exhaustively analyze the toric symmetries of all varieties which arise as toric blowups of $\CC \PP^{3}$. We determine which symmetries act trivially. Furthermore, we explicitly determine the action on homology of each nontrivial toric symmetry.

While toric symmetries and $\CC \PP^{3}$ are of basic interest, our motivation is to study toric symmetry as a computational tool in Gromov--Witten (GW) theory and Donaldson--Thomas (DT) theory. We now attempt to elucidate this idea.

Toric threefolds are a shining example of success in GW Theory and DT theory. Indeed, the all genus
GW/DT theory of any toric threefold $X$ may be computed using multiple techniques. These techniques vary in difficulty and practicality of application, and harness the underlying toric structure of $X$ in different ways.

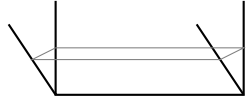
\begin{figure}
\begin{tikzpicture} [thick, scale=1.25]
\draw (0,1) -- (0,0) -- (2,0) -- (2,1);
\draw (-.5,.75) -- (0,0) (1.5,.75) -- (2,0);
\draw[thin,gray] (0,1/2) -- (-1/4,3/8) -- (1.75,3/8) -- (2,.5) -- cycle; 
\end{tikzpicture}
\caption{\small{The polytope of local $\PP^{1}$}}
\label{fig: local-p1-polytope}
\end{figure}

First, any toric threefold $X$ admits analysis via (virtual) fixed point localization. This idea goes back to Kontsevich~\cite{Ko95} and localization of virtual classes was proved by Graber-Pandharipande~\cite{GP}. Localization reduces the defining integrals in GW theory to finite sums. The input data includes the torus fixed points in $X$, which are encoded in the fan (or polytope) of $X$; see Figure~\ref{fig: local-p1-polytope}. While this approach theoretically determines the GW theory of $X$, it is intractable for many choices of toric threefold in practice.

\begin{figure}[h]
\begin{tikzpicture}[thick, scale=.65]
\draw  (0,-1.5) -- (0,0) -- (-1.5,0); 
\draw (0,0) -- (2,2); 
\draw  (3.5,2) -- (2,2) -- (2,3.5);
\end{tikzpicture}
\caption{\small{The web diagram of local $\PP^{1}$}}
\label{fig: local-p1-web}
\end{figure}
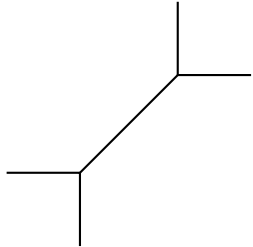

With this in mind, Aganagic-Klemm-Mari\~no-Vafa developed in~\cite{AKMV} the {\emph{topological vertex}}, a tour de force efficiently expressing the GW generating function of a toric Calabi-Yau threefold (CY3) $X$ in terms of combinatorial data encoded in the toric structure of $X$. The combinatorial data in question consists of the (decorated) toric web diagram of $X$; see Figure~\ref{fig: local-p1-web}.

The topological vertex of~\cite{AKMV} is physical, in that it relies on Chern--Simons/Topological String Large-$N$ duality, which remains to be understood mathematically. To overcome this and other difficulties, Li-Liu-Liu-Zhou developed the celebrated {\emph{mathematical theory}} of the topological vertex~\cite{LLLZ}.

The mathematical and physical topological vertices were shown to correspond for most cases already in~\cite{LLLZ}, and therein were conjectured to hold for all toric CY3's. The topological vertex conjecture was shown to hold for an infinite family of toric  local CY3's by Karp-Liu-Mari\~no~\cite{KLM}, and the general proof was given by Maulik-Oblomkov-Okounkov-Pandharipande in~\cite{MOOP}.

Analogous to the topological vertex of~\cite{LLLZ} which computes the GW invariants of any toric CY3, Maulik-Nekrasov-Okounkov-Pandharipande in~\cite{MNOP2} formulate a DT vertex, which they use to establish GW/DT duality for toric CY3's. The trivalent vertex of~\cite{MOOP} computes the GW/DT theory of any toric threefold, proving GW/DT duality for all toric threefolds. 

\begin{figure}[h]
\begin{tikzpicture}[thick, scale=.65]
\draw  (0,-1.5) -- (0,0) -- (-1.5,0); 
\draw (0,0) -- (2,2); 
\draw  (3.5,2) -- (2,2) -- (2,3.5);
\draw (-1.5,0) ellipse (.25cm and .5cm);
\draw (0,-1.5) ellipse (.5cm and .25cm);
\draw (2,3.5) ellipse (.5cm and .25cm);
\draw (3.5,2) ellipse (.25cm and .5cm);
\draw (-1.5,.5) .. controls (0,0) and (2,2) .. (1.5,3.5);
\draw (-1.5,-.5) .. controls (0,0) .. (-.5,-1.5);
\draw (.5,-1.5) .. controls (0,0) and (2,2) .. (3.5,1.5);
\draw (2.5,3.5) .. controls (2,2) .. (3.5,2.5);
\end{tikzpicture}
\caption{\small{The fattened web diagram of local $\PP^{1}$}}
\label{fig: local-p1-fatweb}
\end{figure}
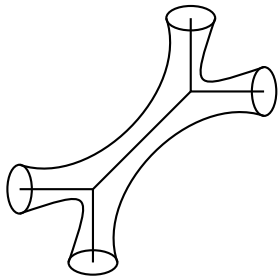

The remodeling conjecture of Bouchard-Klemm-Mari\~no-Pasquetti~\cite{BKMP} is yet another technique applicable to toric CY3's. The remodeling conjecture produces relations among GW invariants difficult to detect via other techniques and takes as input the fattened toric web diagram of $X$; see Figure~\ref{fig: local-p1-fatweb}.

\begin{figure}[h]
\tdplotsetmaincoords{70}{120}
\begin{tikzpicture}[thick,tdplot_main_coords]
	\shade[inner color=green,outer color=white,opacity=.4] (0,0,0) -- (0,0,1) -- (-1,1,1) -- cycle;
	\shade[inner color=green,outer color=white,opacity=.4] (0,0,0) -- (0,1,0) -- (-1,1,1) -- cycle;
	\shade[inner color=green,outer color=white] (0,0,0) -- (0,1,0) -- (1,0,0) -- cycle;
	\draw[->] (0,0,0) -- (-1,1,1) node [anchor=south] {$(-1,1,1)$}; 
	\shade[inner color=green,outer color=white] (0,0,0) -- (0,0,1) -- (1,0,0) -- cycle;
	\shade[inner color=green,outer color=white,opacity=.7] (0,0,0) -- (0,1,0) -- (0,0,1) -- cycle;
	\draw[->] (0,0,0) -- (0,1,0) node [anchor=west] {$(0,1,0)$};
	\draw[->] (0,0,0) -- (0,0,1) node [anchor=south] {$(0,0,1)$};
	\draw[->] (0,0,0) -- (1,0,0) node [anchor=east] {$(1,0,0)$};
\end{tikzpicture}
\caption{\small{The fan of local $\PP^{1}$}}
\label{fig: local-p1-fan}
\end{figure}
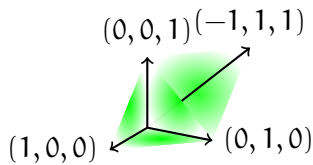

Toric symmetry adds to this list of techniques. Manifest at the level of invariants, as opposed to generating functions, toric symmetry is not detected via remodeling and vertex techniques. The approach of toric symmetry is basic; we identify nontrivial symmetries of the underlying fan or polytope of $X$, and make use of the basic fact that GW/DT invariants are functorial under isomorphism. The fan of our running example, the resolved conifold $\cO(-1) \oplus_{\PP^{1}} \cO (-1)$, is shown in Figure~\ref{fig: local-p1-fan}. 

To illustrate, consider the following example in $\CC \PP^3$.
Let $X$ be the blowup of $\CC \PP^{3}$ at two points $p_{1},p_{2}$ followed by a chain of three lines $l_{1},l_{2},l_{3}$ such that $p_{1}\in l_{1}$ and $p_{2}\in l_{3}$. We call such spaces Class C; for details see Theorem~\ref{thm: classC} below. Let $h$ be the pullback of the line class in $\CC \PP^{3}$, let $e_{i}$ denote the class of a line in the exceptional divisor above $p_{i}$, and let $f_{i}$ denote the fiber class in the exceptional divisor above $l_{i}$. 

What is $GW^{X}_{0,f_{1}-f_{2}}$? It may be nontrivial to compute this invariant using the topological vertex or localization. However, using toric symmetry, in the form of Theorem~\ref{thm: classC} below, we immediately obtain 
\begin{equation*}
\boxed{
GW^{X}_{0,f_{1}-f_{2}} = GW^{X}_{0,h-e_{1}-e_{2}} =1}
\end{equation*}
where the second equality holds as the invariant in question counts the number of lines in space through the two points $p_{1}$ and $p_{2}$.

The example first discovered of toric symmetry in GW theory was \emph{Cremona Symmetry}, which is the symmetry of $\PP^{n}$ induced via (the resolution of) the classical Cremona birational map on $\PP^{n}$. The idea to use the stability of GW invariants (of surfaces) under birational transformations goes back at least to Miranda-Crauder; see the final remarks of Section 11 in~\cite{CrMi}. This philosophy was actualized on $\PP^{2}$ by G\"ottsche-Pandharipande~\cite{GoPa} and also Bryan-Leung~\cite{BL}. The Cremona symmetry of $\PP^{3}$ was studied by Gathmann~\cite{Ga01} and Bryan-Karp~\cite{BK}. 

The present work exhaustively identifies all toric symmetries of toric blowups of $\PP^{3}$, of which Cremona symmetry is a single example.

\subsection*{Main Results}

In slightly more detail: given the fan $\Sigma\subset \ZZ^n$ of a toric variety $X$, recall that a \emph{toric symmetry} of $X$ is an automorphism of the lattice $\ZZ^n$ which permutes the cones of $\Sigma$. Moreover, by {\emph{nontrivial}} toric symmetry we refer to those toric symmetries which act nontrivially on GW/DT invariants/homology; see Section~\ref{sec: toric_symmetries} for further discussion.

\begin{theorem}\label{thm: symmetries}
There exist precisely four classes of toric blowups of $\PP^3$ which have nontrivial toric symmetry. These four classes, labelled A, B, C and D, are described in Theorems~\ref{thm: classA}, \ref{thm: classB}, \ref{thm: classC}, and~\ref{thm: permutohedron} respectively.

Moreover, a space of Class A, B or D admits a unique toric symmetry (up to relabeling), whereas there are precisely four distinct nontrivial symmetries for Class C varieties.
\end{theorem}

\begin{remark}\label{rem: graph_associahedra}
The blowup of $\PP^3$ at its four torus fixed points and subsequently the proper transforms of its six torus invariant lines yields the toric variety $\hat X$ whose polytope $\Delta_{\hat X}$ is the permutohedron $\Pi_3$. The automorphism of $\hat{X}$ induced by reflection of $\Pi_{3}$ through the origin resolves the classical Cremona transform on $\PP^{3}$; see~\cite[Section 6]{Ga01}. An analysis of this reflection yields Theorem~\ref{thm: permutohedron} below and is given in~\cite{BK}.

The permutohedron is one among a class of polytopes known as graph associahedra which include the cyclohedron and associahedron; for the definition and basic properties of graph associahedra, see~\cite{Dev06}. This class of polytope is encompassed by polytopes $\Delta_X$ where $X$ is a toric blowup (in the sense below) of $\PP^3$~\cite{Dev09}. It then follows immediately from Theorem~\ref{thm: symmetries} that the permutohedron is the only graph associahedron which admits a nontrivial toric symmetry.
\end{remark}

We now construct these four families and describe their toric symmetries. In what follows $X$ will be an iterated blowup of $\mathbb{P}^3$ at a specified configuration of points and lines. Throughout this work, we say
the variety $X$ is a {\emph{toric blowup}} of $\PP^{3}$ if $X$ is an iterated blowup of $\PP^{3}$ only along torus invariant subvarieties of $\PP^{3}$ (or their proper transforms). In particular, we are not interested in spaces obtained by blowups with centers in the exceptional locus, as their geometry is far from $\PP^{3}$.

 Abusively, we often identify divisors and their classes; when we need care, classes will be decorated with brackets, e.g. $[D]$. We denote by $H$ the pullback to $X$ of the hyperplane class in $\mathbb{CP}^3$. We denote by $E_\alpha$ the exceptional divisor above a point $p_\alpha$ and by $F_{\alpha'}$ the exceptional divisor above a line $\ell_{\alpha'}$ for appropriate indices $\alpha, \alpha'$. Further, we will let $h$ and $e_\alpha$ denote the classes of a line in the divisors $H$ and $E_\alpha$ respectively, and $f_{\alpha'}$ denotes the fibre class in $F_{\alpha'}$. The homology groups $H_4(X;\ZZ)$ and $H_2(X;\ZZ)$ are spanned by divisor and curve classes, respectively:
\[
H_4(X;\ZZ) = \left\langle H, E_\alpha, F_{\alpha'} \right\rangle, \ \ H_2(X;\ZZ) = \left\langle h, e_\alpha, f_{\alpha'} \right\rangle.
\]
For computation of the chow ring, see Section~\ref{sec: chow_ring} below. 

\begin{theorem}\label{thm: classA}
Let $X$ be the blowup of $\mathbb{P}^3$ along a  point $p$ and two intersecting distinct lines $\ell_1$ and $\ell_2$, such that $p \neq \ell_1 \cap  \ell_2$. We call such a space a \emph{Class A} blowup; see Figure~\ref{fig: class-A-blowup-diagram}. Then, given $\beta = dh - a_1e - a_2f_1-a_3f_2 \in H_{2}(X;\ZZ)$, there exists a unique nontrivial toric symmetry $\tau_{A}$ on $X$, and its action on homology is given by
\[
(\tau_{A})_{*} \beta = \beta '
\]
where $\beta' = d'h - a_1'e - a_2'f_1-a_3'f_2$ has coefficients given by
\begin{eqnarray*}
d' &=& 2d-a_1-a_2-a_3\\
a_1' &=& d-a_2-a_3\\
a_2' &=& d-a_1-a_3\\
a_3' &=& d-a_1-a_2. 
\end{eqnarray*}

\begin{figure}[h!]
\begin{tikzpicture}[thick]
	\draw (1,0) -- (2,1) node [midway,anchor=east] {$\ell_1$};
	\draw (2,1) -- (3,0) node [midway,anchor=west] {$\ell_2$};
	\filldraw (2,-.1) circle (1pt) node [anchor=north] {$p$};
\end{tikzpicture}
\caption{\small{The ordered Class A blowup center}}
\label{fig: class-A-blowup-diagram}
\end{figure}
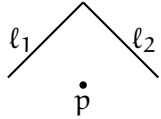

\end{theorem}

\begin{theorem}\label{thm: classB}
Let $X$ be the sequential blowup of $\mathbb{P}^3$ at distinct points $p_1$ and $p_2$ and three pairwise intersecting lines $\ell_1,\ell_2$ and $\ell_3$ such that $p_1\in \ell_1$, $p_2\in \ell_3$ and $p_1,p_2\notin \ell_2$. Then $X$ is called a \emph{Class B} blowup; see Figure~\ref{fig: class-B-blowup-diagram}.  Then, given $\beta = dh - a_1e_1 - a_2e_2-a_3f_1-a_4f_2-a_5f_3$, there exists a unique toric symmetry $\tau_{B}$ on $X$ and its action on homology is given by
\[
(\tau_{B})_{*}\beta =\beta '
\]
where $\beta' = d'h - a_1'e_1 - a_2'e_2-a_3'f_1-a_4'f_2-a_5'f_3$ has coefficients given by
\begin{eqnarray*}
d' &=& 2d-a_1-a_2-a_3-a_4\\
a_1' &=& a_5\\
a_2' &=& d-a_1-a_3-a_4\\
a_3' &=& d-a_2-a_4-a_5\\
a_4' &=& d-a_1-a_2-a_3\\
a_5' &=& a_1.
\end{eqnarray*}

\begin{figure}[h!]
\begin{tikzpicture}[thick]
	\filldraw (-.5,0) circle (1pt) node [anchor=north] {$p_1$};
	\filldraw (.5,0) circle (1pt) node [anchor=north] {$p_2$};
	\draw (-.5,0) -- (-1,1) node [midway,anchor=east] {$\ell_1$};
	\draw (-1,1) -- (1,1) node [midway,anchor=south] {$\ell_2$};
	\draw (1,1) -- (.5,0) node [midway,anchor=west] {$\ell_3$};
\end{tikzpicture}
\caption{\small{The ordered Class B blowup center}}
\label{fig: class-B-blowup-diagram}
\end{figure}
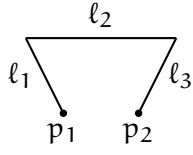

\end{theorem}

\begin{theorem}\label{thm: classC}
Let $X$ be the sequential blowup of $\mathbb{P}^3$ at distinct points $p_1$ and $p_2$ and three pairwise intersecting lines $\ell_1,\ell_2$ and $\ell_3$ such that $p_1\in \ell_2$, $p_2\in \ell_3$ and $p_1,p_2\notin \ell_1$. We term such spaces \emph{Class C} blowups; see Figure~\ref{fig: class-C-blowup-diagram}. Then, there exist precisely four nontrivial toric symmetries $\tau_{C}, \sigma_{C}, \sigma^{2}_{C},$ and $\sigma_{C} \tau_{C}$ of $X$. With $\beta = dh - a_1e_1 - a_2e_2-a_3f_1-a_4f_2-a_5f_3$, their action on cohomology is given by 
\[
(\tau_{C})_{*}\beta =\beta' \qquad\qquad (\sigma_{C})_{*}\beta = \beta''
\]
where $\beta' = d'h - a_1'e_1 - a_2'e_2-a_3'f_1-a_4'f_2-a_5'f_3$ has coefficients given by
\begin{eqnarray*}
d' &=& 2d-a_1-a_2-a_3-a_4\\
a_1' &=& a_5\\
a_2' &=& d-a_2-a_3-a_4\\
a_3' &=& d-a_1-a_2-a_4\\
a_4' &=& a_2\\
a_5' &=& d-a_1-a_3-a_5,
\end{eqnarray*}
and $\beta'' = d''h - a_1''e_1 - a_2''e_2-a_3''f_1-a_4''f_2-a_5''f_3$ has coefficients given by
\begin{eqnarray*}
d'' &=& 2d-a_1-a_2-a_3-a_5\\
a_1'' &=& a_4\\
a_2'' &=& d-a_1-a_3-a_5\\
a_3'' &=& d-a_1-a_2-a_5\\
a_4'' &=& a_1\\
a_5'' &=& d-a_2-a_3-a_4.
\end{eqnarray*}

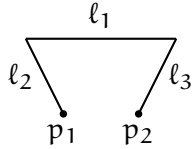
\begin{figure}[h!]
\begin{tikzpicture}[thick]
	\filldraw (-.5,0) circle (1pt) node [anchor=north] {$p_1$};
	\filldraw (.5,0) circle (1pt) node [anchor=north] {$p_2$};
	\draw (-.5,0) -- (-1,1) node [midway,anchor=east] {$\ell_2$};
	\draw (-1,1) -- (1,1) node [midway,anchor=south] {$\ell_1$};
	\draw (1,1) -- (.5,0) node [midway,anchor=west] {$\ell_3$};
\end{tikzpicture}
\caption{\small{The ordered Class C blowup center}}
\label{fig: class-C-blowup-diagram}
\end{figure}
\end{theorem}

\begin{theorem}[Bryan-Karp~\cite{BK}]\label{thm: permutohedron}
Let $X$ be the sequential blowup of $\PP^3$ at 4 distinct points $p_1,\ldots,p_4$ and the six distinct lines $\ell_{ij}$ between them. Let $\beta$ be  given by
\[
 \beta = dh - \sum_{i=1}^{4} a_ie_i - \sum_{1\leq i< j\leq 4} b_{ij}f_{ij} \in H_{2}(X;\ZZ).
\]
There exists a unique toric symmetry $\tau_{D}$ of $X$, and its action on homology is given by $(\tau_{D})_{*} \beta =\beta '$,
where $\beta' = d'h - \sum_{i} a_i'e_i - \sum_{ij} b_{ij}'f_{ij}$ has coefficients given by 
\begin{eqnarray*}
d' &=& 3d - 2\sum_{i=1}^{4} a_i \\
a_i' &=& d - a_j - a_k - a_l-b_{ij}-b_{ik}-b_{il}\\
b_{ij}' &=& b_{kl},
\end{eqnarray*}
where $\{i,j,k,l\} = \{1,2,3,4\}$. We refer to $X$ as a \emph{Class D} blowup of $\PP^{3}$; see Figure~\ref{fig: class-D-blowup-diagram}. 

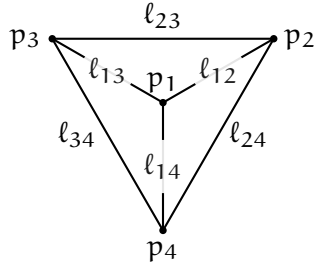
\begin{figure}[h!]
\def\lineseg{1.7}
\begin{tikzpicture}[thick]
	\draw (30:\lineseg) -- (150:\lineseg) node [midway,anchor=south] {$\ell_{23}$};
	\draw (150:\lineseg) -- (270:\lineseg) node [midway,anchor=east] {$\ell_{34}$};
	\draw (270:\lineseg) -- (30:\lineseg) node [midway,anchor=west] {$\ell_{24}$};

	\draw (30:\lineseg) -- (0,0) node [midway,shape=circle,inner sep=1,fill=white,opacity=0.9,text opacity=1] {$\ell_{12}$};
	\draw (150:\lineseg) -- (0,0) node [midway,shape=circle,inner sep=1,fill=white,opacity=0.9,text opacity=1] {$\ell_{13}$};
	\draw (270:\lineseg) -- (0,0) node [midway,shape=circle,inner sep=1,fill=white,opacity=0.9,text opacity=1] {$\ell_{14}$};

	\filldraw (0,0) circle (1pt) node [anchor=south] {$p_1$};
	\filldraw (30:\lineseg) circle (1pt) node [anchor=west] {$p_2$};
	\filldraw (150:\lineseg) circle (1pt) node [anchor=east] {$p_3$};
	\filldraw (270:\lineseg) circle (1pt) node [anchor=north] {$p_4$};
\end{tikzpicture}
\caption{\small{The ordered Class D blowup center}}
\label{fig: class-D-blowup-diagram}
\end{figure}

\end{theorem}

\begin{corollary}\label{main corollary}
Let $X$ be a toric blowup of $\CC \PP^{3}$ of Class $A$, $B$, $C$, or $D$. Then we have the following equality of GW and DT invariants
\begin{align}\label{eq: main corollary}
GW^{X}_{g,\beta} &= GW^{X}_{g,\beta '}\\
DT^{X}_{g,\beta} &= DT^{X}_{g,\beta '},
\end{align}
where $\beta$ and $\beta '$ (or even $\beta''$) are given in Theorems~\ref{thm: classA} through~\ref{thm: permutohedron} above.
\end{corollary}

\begin{remark}\label{rem: higher-virtual-dimension}
The results of Corollary~\ref{main corollary} easily extend to invariants with insertions (the case of higher virtual dimension); see Section~\ref{subsec: nonzero_vd}.
\end{remark}

A lower dimensional analog of Theorem~\ref{thm: symmetries} is given below, and is proved using analogous methods. 

\begin{theorem}\label{thm: symmetriesP2}
There exists precisely one family of spaces, obtained from the maximal toric blowup of $\PP^2$, which have toric symmetries that induce nontrivial action on Chow, and lift to nontrivial relations of GW invariants. In particular, this is the blowup of $\PP^2$ at its torus fixed points. The unique toric symmetry is the Cremona Symmetry studied in G\"ottsche-Pandharipande~\cite{GoPa}. 
\end{theorem}

\subsection*{Acknowledgments}
We gratefully acknowledge Jim Bryan and Renzo Cavalieri for useful conversations and comments on an earlier draft. This work was partially supported by National Science Foundation Grant DMS-083996 and the Beckman Foundation.

\section{GW Theory and DT Theory}\label{sec: GW/DT}
We now recall Gromov-Witten theory, Donaldson-Thomas theory, and fix notation. 

\subsection{GW Theory}
Let $X$ be a smooth projective threefold. For a given curve class $\beta \in H_{2}(X;\ZZ )$, the moduli stack of genus-$g$, $n$-pointed stable maps to $X$ is denoted $\Mbar_{g,n }(X,\beta )$. Elements of $\Mbar_{g,n}(X,\beta)$ are isomorphism classes of maps 
\[
f: (C, p_{1},\dotsc ,p_{n}) \rightarrow X,
\]
where $C$ is a genus $g$ curve with at worst nodal singularities and $p_{i}\in C$ are marked points. Such a map is stable if its automorphism group is finite, where an isomorphism of the map $f: (C, p_{1},\dotsc ,p_{n}) \rightarrow X$ is given by a commutative triangle. 
\[
\xymatrix{ (C,p_{1},\dotsc ,p_{n}) \ar[rr]^{\tau} \ar[dr]_{f} & &  (C',p'_{1},\dotsc,p'_{n})  \ar[dl]^{f'} \\
           & X &  }
\]
Here, $\tau:C \rightarrow C'$ is an isomorphism of curves respecting marked points, $\tau (p_i)=p'_{i}$.

The moduli stack $\Mbar_{g,n}(X,\beta )$ has a virtual fundamental class
\[
[\Mbar_{g,n}(X,\beta )]^{vir} \in H_{vd}(\Mbar_{g,n}(X,\beta );\QQ),
\]
where the virtual dimension $vd$ is given by 
\begin{align*}
vd &= (\dim X -3)(1-g) -K_{X}\cdot \beta +n\\
&=n-K_{X}\cdot \beta.
\end{align*}
In particular, the moduli stack of unmarked stable maps to a CY3 has a virtual fundamental zero-cycle. For construction of the virtual class see~\cite{Be97,BeFa97}.

Given a collection of cohomology classes $\gamma_{1},\dotsc ,\gamma_{n} \in H^{*}(X)$, GW invariants are constructed to virtually count cycles in $X$ of class $\beta$ intersecting each Poincar\'e dual $[\gamma_{i}]^{PD}$. This is achieved by pulling back the classes $\gamma_{i}$ and pairing with the virtual class. We pullback via evaluation morphisms $\ev_{i}$ given by
\begin{align*}
\ev_{i}: \Mbar_{g,n}(X,\beta ) &\longrightarrow X\\
[f:(C ,p_{1},\dotsc ,p_{n})\rightarrow X] &\longmapsto f(p_{i}).
\end{align*}

The genus-$g$ class $\beta$ GW invariant of $X$ with insertions $\gamma_{1},\dotsc ,\gamma_{n}$ is defined by
\[
GW(\gamma_{1},\dotsc ,\gamma_{n})^{X}_{g,\beta} = \int_{[\Mbar_{g,n}(X,\beta )]^{vir}} \prod_{i=1}^{n} \ev_{i}^{*} \gamma_{i}. 
\]
GW invariants are invariant under symplectic deformation of $X$. In particular, they are well defined at the level of isomorphism class. More precisely, if $\tau : X\longrightarrow X$ is an automorphism, then
\begin{equation}\label{eq: gw-symmetry}
GW(\gamma_{1},\dotsc ,\gamma_{n})^{X}_{g,\tau_{*}\beta } = GW( \tau^{*}\gamma_{1},\dotsc ,\tau^{*}\gamma_{n} )^{X}_{g,\beta }.
\end{equation}
For further details concerning the definition of GW invariants, see the text by Hori et al~\cite{Hori-et-al}.

\subsection{DT Theory} Donaldson-Thomas invariants also virtually count curves in a smooth projective threefold $X$ of class $\beta \in H_{2}(X;\ZZ )$ intersecting Poincar\'e duals of cohomology classes $\gamma_{1},\dotsc ,\gamma_{r} \in H^{*}(X)$. Instead of probing the geometry of $X$ via stable maps as in GW theory, rather the basic idea in DT theory is to study sheaves on $X$.

For an ideal sheaf\footnote{We follow the conventions of~\cite{MNOP2,MOOP}. An ideal sheaf is a rank 1 torsion free sheaf with trivial determinant.} $\cI$, there exists an injection into its double dual
\[
0 \longrightarrow \cI  \longrightarrow \cI ^{\vee \vee}.
\]
But 
\[
\cI^{\vee \vee} \cong \cO_{X},
\]
so $\cI$ determines a subscheme $Y$ given by
\[
0 \longrightarrow \cI \longrightarrow \cO_{X} \longrightarrow \cO_{Y} \longrightarrow 0.
\]
Since $\cI$ has trivial determinant, $Y$ has components of dimension zero and one. The weighted one dimensional components of $Y$ determine a homology class 
\[
[Y] \in H_{2}(X;\ZZ  ).
\]

The moduli space of ideal sheaves $\cI$ with holomorphic Euler characteristic $\chi (\cO_{Y})=n$ and class $[Y]=\beta \in H_{2}(X;\ZZ )$ is denoted $I_{n}(X,\beta )$. Similar to GW invariants, DT invariants are defined by integrating against the virtual class $[I_{n}(X,\beta)]^{\vir }$ of dimension
\[
\dim [I_{n}(X,\beta )]^{\vir} = \int_{\beta}c_{1} (T_{X}). 
\]
The construction of this virtual class and other foundational aspects of DT theory may be found in~\cite{MNOP1,Thomas}. 

In order to integrate against the virtual class, we need to pull back the classes $\gamma_{i}$ from $X$ to $I_{n}(X,\beta )$. This is done using the universal ideal sheaf and the associated universal subscheme. 

By~\cite[Section 2.2]{MNOP2} and~\cite[Section 1.2]{MOOP} and there exists a universal ideal sheaf
\[
\fI \longrightarrow I_{n}(X,\beta )\times X
\]
with well defined Chern classes\footnote{The second Chern class $c_{2}(\fI )$ is interpreted as the universal subscheme~\cite{MOOP}.}. Let $\pi_{i}$ denote the respective projection maps. The DT invariants are defined by push-pulling Chern classes via $\pi_{i}$. For each $\gamma \in H^{*}(X)$, define the operator $c_{2}(\gamma )$ by, for any $\xi \in H_{*}(I_{n}(X,\beta ))$,
\begin{equation}\label{eq: dt-pullback}
c_{2}(\gamma )(\xi) = \pi_{1*} \left( c_{2}(\fI )\cdot \pi_{2}^{*}(\gamma )\cap \pi_{1}^{*}(\xi) \right).
\end{equation}
For details of this construction, including the pullback of the homology class $\xi$ in Equation~\ref{eq: dt-pullback}, see~\cite[Section 1.2]{MOOP}. 

The class $(n,\beta )$ DT invariant of $X$ with insertions $\gamma_{1},\dotsc \gamma_{r}$ is defined by
\[
DT^{X}_{n,\beta}(\gamma_{1},\dotsc ,\gamma_{r}) = \langle \gamma_{1},\dotsc ,\gamma_{r}\rangle^{X}_{n,\beta} = \int_{[I_{n}(X,\beta )]^{\vir}} \prod_{i=1}^{r} c_{2}(\gamma_{i}). 
\]
The DT invariants of $X$ are indeed deformation invariants and in particular are well defined at the level of automorphism. So, for any automorphism $\tau$ of $X$, we have 
\begin{equation}\label{eq: dt-symmetry}
DT^{X}_{n,\tau_{*}\beta}(\gamma_{1},\dotsc ,\gamma_{r}) =
DT^{X}_{n,\beta}(\tau^{*} \gamma_{1},\dotsc ,\tau^{*}\gamma_{r}).
\end{equation}

\section{Toric Blowups of $\CC \PP^{3}$}\label{sec: toric_blowups}
In this section, we explicitly construct toric models of the varieties in question, and as a byproduct obtain a presentation of their (co)homology.

\subsection{Iterated Toric Blowups}\label{sec: iterated_toric_blowups}
Recall that the fan $\Sigma_{\PP^3}\subset \ZZ^3$ of $\PP^3$ has one-skeleton with primitive generators
\[
\left.
\begin{array}{l@{\qquad}c}
v_1 = (-1,-1,-1),  &  v_2 = (1,0,0),\\
 v_3 = (0,1,0), & v_4 = (0,0,1).
\end{array}
\right.
\]
and maximal cones given by 
\[
\left.
\begin{array}{cc}
\langle v_1,v_2,v_3 \rangle,  &  \langle v_1,v_2,v_4 \rangle,\\
\langle v_1,v_3,v_4 \rangle, & \langle v_2,v_3,v_4 \rangle.
\end{array}
\right.
\]
The lower dimensional cones of $\Sigma_{\PP^3}$ can be found by intersecting maximal cones. Also note that the fan $\Sigma_{\PP^3}$ is the normal fan over the faces of a polytope $\Delta_{\PP^3}$, which has the adjacency structure of a 3-simplex. 

\subsubsection{Notation}\label{sec: tfixed_notation} We will label torus invariant subvarieties by the primitive generators of their cones. Specifically, the point $p_{ijk}$ will refer to the orbit closure of the cone $\langle v_i,v_j,v_k \rangle$. Similarly, the line $\ell_{ij}$ will refer to the orbit closure of the cone $\langle v_i,v_j \rangle$. We will label a new element of the 1-skeleton, introduced to subdivide the cone $\sigma = \langle v_i,\ldots,v_j \rangle$, by $v_{i \cdots j}$.

\begin{remark}\label{rem: proper-transforms}
We are blowing up proper transforms of subvarieties sequentially, and thus, any toric blowup of $\PP^3$ is determined by a collection of $T$-fixed points and an {\it ordered} collection of $T$-fixed lines in $\PP^3$, chosen from the set $\{p_{ijk},\ell_{rs}\}$ for $\{i,j,k,r,s\} = \{1,2,3,4\}$. 
\end{remark}

\subsubsection{Fans}\label{sec: blowup_configs} With the above notation, we will now explicitly construct the fan of a toric blowup of $\PP^3$ of Class A. The subvarieties we will blow up are $p_{123}$, $\ell_{34}$ and $\ell_{24}$. To blow up $\PP^3$ at $p_{123}$, we subdivide the cone $\langle v_1,v_2,v_3 \rangle$ by inserting a ray generated by $v_{123} = v_1+v_2+v_3 = (0,0,-1)$, and replace the single maximal cone with three cones, according to the subdivision
\[
\langle v_1,v_2,v_3 \rangle \to \langle v_1,v_2,v_{123} \rangle, \ \langle v_1,v_3,v_{123} \rangle, \ \langle v_2,v_3,v_{123} \rangle.
\]
This new fan, $\Sigma_1$, determines the toric variety of $\PP^3$ blown up at $p_{123}$. The proper transform of $\ell_{34}$ corresponds to the cone $\langle v_3,v_4 \rangle$. So, to blow up the line $\ell_{34}$ we subdivide $\langle v_3,v_4 \rangle$, we insert the ray generated by $v_{34} = v_3+v_4 = (0,1,1)$. We point out that the subdivision of any sub maximal cone, in this case $\langle v_3,v_4 \rangle$, necessitates the subdivision of any cones that may contain it as a proper face. In this case, 
\begin{eqnarray*}
\langle v_3,v_4 \rangle &\to& \langle v_3,v_{34} \rangle, \ \langle v_4,v_{34} \rangle. \\
\langle v_1,v_3,v_4 \rangle &\to& \langle v_1,v_3,v_{34} \rangle, \ \langle v_1,v_4,v_{34} \rangle \\
\langle v_2,v_3,v_4 \rangle &\to& \langle v_2,v_3,v_{34} \rangle, \ \langle v_2,v_4,v_{34} \rangle. 
\end{eqnarray*}
The above process yields a new fan $\Sigma_2$, whose toric variety $X_{\Sigma_2}$ is the blowup of $X_{\Sigma_1}$ at $\ell_{34}$. We now repeat the process with the proper transform of the line $\ell_{24}$. The process yields a fan $\hat\Sigma$, whose toric variety is the blowup of $X_{\Sigma_2}$ at $\ell_{24}$, or in other words the toric blowup of $\PP^3$ at $p_{123}$, $\ell_{34}$ and $\ell_{24}$. 

\begin{figure}
\small
\tdplotsetmaincoords{59}{118}
\newsavebox{\tempbox}
\sbox{\tempbox}{
	\begin{tikzpicture}[thick,tdplot_main_coords]
	    \draw[->] (0,0,0) -- (1,0,0) node[anchor=north east]{$v_2,\ (1,0,0)$};
	    \draw[->] (0,0,0) -- (0,1,0) node[anchor=west]{$(0,1,0),\ v_3$};
	    \draw[->] (0,0,0) -- (0,0,1) node[anchor=south]{\makebox[0pt]{$v_4, \hspace{15pt}$} $(0,0,1)$};
	    \draw[->] (0,0,0) -- (-1,-1,-1) node[anchor=east]{$v_1,\ (-1,-1,-1)$};
	\end{tikzpicture}
}
\newlength{\tempboxlength}
\settowidth{\tempboxlength}{\usebox{\tempbox}}
\subfloat[(i) One-skeleton of $\PP^3$]{
	\usebox{\tempbox}
}
\makebox[\tempboxlength]{
\hfill
\subfloat[(ii) \ldots blown up at $p_{123}$]{
	\begin{tikzpicture}[thick,tdplot_main_coords]
	    \draw[->] (0,0,0) -- (1,0,0) node[anchor=north east]{$(1,0,0)$};
	    \draw[->] (0,0,0) -- (0,1,0) node[anchor=west]{$(0,1,0)$};
	    \draw[->] (0,0,0) -- (0,0,1) node[anchor=south]{$(0,0,1)$};
	    \draw[->] (0,0,0) -- (-1,-1,-1) node[anchor=east]{$(-1,-1,-1)$};
	    \draw[->] (0,0,0) -- (0,0,-1) node[anchor=north]{$(0,0,-1)$\makebox[0pt]{$\hspace{27pt},\ v_{123}$}};
	\end{tikzpicture}
}
\hfill
}
\subfloat[(iii) \ldots then at $\ell_{34}$]{
	\begin{tikzpicture}[thick,tdplot_main_coords]
	    \draw[->] (0,0,0) -- (1,0,0) node[anchor=north east]{$(1,0,0)$};
	    \draw[->] (0,0,0) -- (0,1,0) node[anchor=west]{$(0,1,0)$};
	    \draw[->] (0,0,0) -- (0,0,1) node[anchor=south]{$(0,0,1)$};
	    \draw[->] (0,0,0) -- (-1,-1,-1) node[anchor=east]{$(-1,-1,-1)$};
	    \draw[->] (0,0,0) -- (0,0,-1) node[anchor=north]{$(0,0,-1)$};
	    \draw[->] (0,0,0) -- (0,1,1) node[anchor=west]{$(0,1,1),\ v_{34}$};
	\end{tikzpicture}
}
\hfill
\makebox[\tempboxlength]{
\hfill
\subfloat[(iv) \ldots then at $\ell_{24}$]{
	\begin{tikzpicture}[thick,tdplot_main_coords]
	    \draw[->] (0,0,0) -- (1,0,0) node[anchor=north east]{$(1,0,0)$};
	    \draw[->] (0,0,0) -- (0,1,0) node[anchor=west]{$(0,1,0)$};
	    \draw[->] (0,0,0) -- (0,0,1) node[anchor=south]{$(0,0,1)$};
	    \draw[->] (0,0,0) -- (-1,-1,-1) node[anchor=east]{$(-1,-1,-1)$};
	    \draw[->] (0,0,0) -- (0,0,-1) node[anchor=north]{$(0,0,-1)$};
	    \draw[->] (0,0,0) -- (0,1,1) node[anchor=west]{$(0,1,1)$};
	    \draw[->] (0,0,0) -- (1,0,1) node[anchor=east]{$v_{24},\ (1,0,1)$};
	\end{tikzpicture}
}
}
\caption{\small{Constructing the fan for the iterated blowup of Theorem \ref{thm: classA}}}
\label{fig: blowup-construction}
\end{figure}
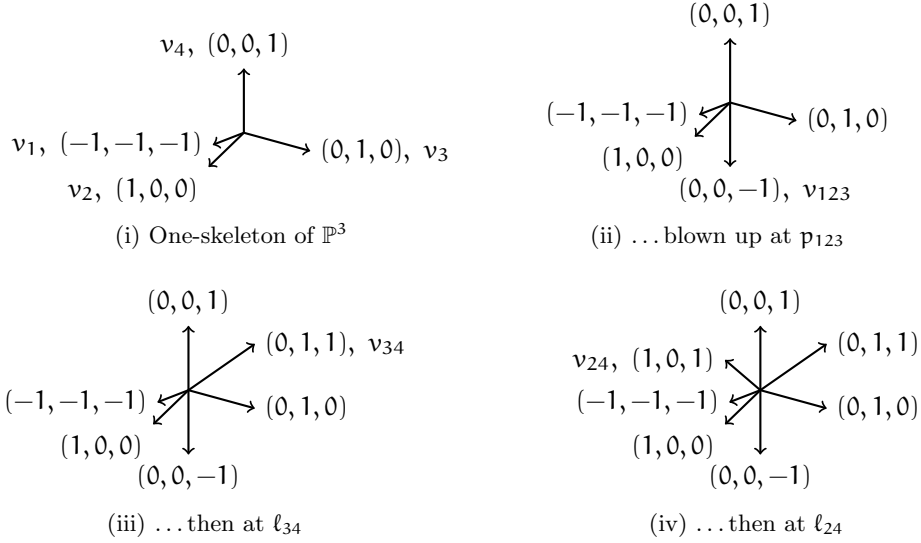

For other spaces, the procedure is above is carried out analogously. Given an ordered collection of torus invariant subvarieties, we subdivide their cones in the above manner to obtain a new fan $\hat \Sigma$ whose toric variety $X_{\hat \Sigma}$ is the blowup of $\PP^3$ at this collection of subvarieties.

\subsection{Chow Ring}\label{sec: chow_ring}
For a nonsingular projective variety $X$ with fan $\Sigma_X$, the Chow ring $A^*(X)$ is a quotient of the polynomial algebra over $\ZZ$ generated by the variables $\{D_\alpha\}$ in bijection with $v_\alpha \in \Sigma^{(1)}$, see \cite{Fu}. We will now explicitly construct the Chow ring for the Class A blowup. Let $X$ be the Class A toric blowup constructed in Section~\ref{sec: blowup_configs}. Then the Chow ring is given by \[
A^*(X) =\ZZ[D_1,D_2,D_3,D_4,D_{123},D_{34},D_{24}]/I.
\]
The ideal $I$ contains the relations 
\begin{enumerate}
\item[(i)] $D_{\alpha_1}\cdots D_{\alpha_k}$ for $v_{\alpha_1},\cdots ,v_{\alpha_k}$ not in a cone of $\Sigma_X$.
\item[(ii)] $\sum_\alpha \langle u,v_\alpha \rangle D_\alpha$ for $u\in \ZZ^n$, a basis element of the lattice. 
\end{enumerate}
The relations of the first type are
\[
D_2D_4, \ D_3D_4, \ D_3D_{24}, \ D_4D_{123}, \ D_{123}D_{34}, \ D_{123}D_{24},
\]
\[
D_1D_2D_3, \ D_1D_2D_{34}, \ D_1D_{34}D_{24}
\]
The relations of the second type are
\[
-D_1+D_2+D_{24}, \ -D_1+D_3+D_{34}, \ -D_1+D_4-D_{123}+D_{34}+D_{24}.
\]

In particular, the Chow ring is generated by divisor classes of the orbit closures of elements in the 1-skeleton. Another, more geometric presentation for $A^*(X)$ is given as follows. For $X$, a toric blowup of $\PP^3$ at points $p_\alpha$ and lines $\ell_{\alpha'}$, let $H$ be the pullback of the hyperplane class in $\PP^3$ to $X$ and $E_\alpha$ and $F_{\alpha'}$ be the exceptional divisors above $p_\alpha$ and $\ell_{\alpha'}$ respectively. Then $A^1(X)$ is generated by the classes $H$, $E_\alpha$ and $F_{\alpha'}$.

We remark that the anti-canonical class is given by
\[
-K_{X} = \sum_\alpha D_\alpha = 4H - 2\sum E_{\alpha} - \sum F_{\alpha'}.
\]

The Chow rings for Classes B, C and D are easily and analogously determined by their fans.

\section{Toric Symmetry}\label{sec: toric_symmetries}

\subsection{Toric Symmetries}\label{subsec: nontrivial}
Recall that a toric symmetry of $X$ is an automorphism $\tau$ of the lattice $\ZZ^n$ which permutes the cones of $\Sigma$. In particular, these automorphisms permute the 1-dimensional cones, $\Sigma^{(1)}$. Since the elements of $\Sigma^{(1)}$ are in bijection with the divisor classes, and further, these classes generate Chow, the unique linear extension of $\tau$ defined on $A^1(X)$ which commutes with the product structure of $A^*(X)$ is the pullback $\tau^*$.

 We are interested in characterizing {\it nontrivial} toric symmetries at the level of the fan. There are three distinct types of divisors classes in a toric blowup of $\PP^3$. The first are those involving the pullback of the hyperplane class. The second and third are the classes of exceptional divisors above the centers of blowups of points and lines respectively. Any symmetry that is nontrivial at the level of GW/DT theory will exchange these three classes of divisors. It can be verified easily that in a toric blowup of $\PP^3$, any nontrivial symmetry maps the class $H$ to a divisor class with coefficient in $H$ greater than 1. 

\subsection{Proof of Theorem~\ref{thm: symmetries}}\label{sec: proof_5} Let $X$ be a toric blowup of $\PP^3$. Let $\Sigma_X$ be the fan of $X$. Notice that the primitive generators of $\Sigma_X$ include the standard basis of $\ZZ^3$. Also observe that the elements of the one-skeleton are sums of $v_1,\ldots,v_4$, as constructed in Section~\ref{sec: toric_blowups}. Since $\tau$ acts on $\Sigma_X^{(1)}$, the standard basis elements of $\ZZ^3$ must be mapped by $\tau$ to elements whose entries are in the set $\{-1,0,1\}$. Therefore the automorphism group of the fan is finite. Moreover, this analysis yields a computational method to determine this group of lattice isomorphisms. From Section~\ref{sec: chow_ring}, we see that the action of $\tau$ on $\Sigma_X^{(1)}$ yields a map $\tau^*$ on $A^*(X)$.

The sets $\{D_1,\ldots, D_4\}$, $\{D_{ijk}\}$ and $\{D_{rs}\}$ correspond to $H$'s, $E_{\alpha}$'s and $F_{\alpha'}$'s, respectively. Thus, nontrivial toric symmetries are characterized by those which exchange elements of these three sets above amongst each other. 

Since the number of toric blowups is finite, this characterization for maps which pushforward nontrivially to $A_*(X)$ allows us to computationally find all nontrivial toric symmetries. We have written Sage code and actualized this strategy. The code itself is contained in the source of this document.

\begin{figure}[h!]
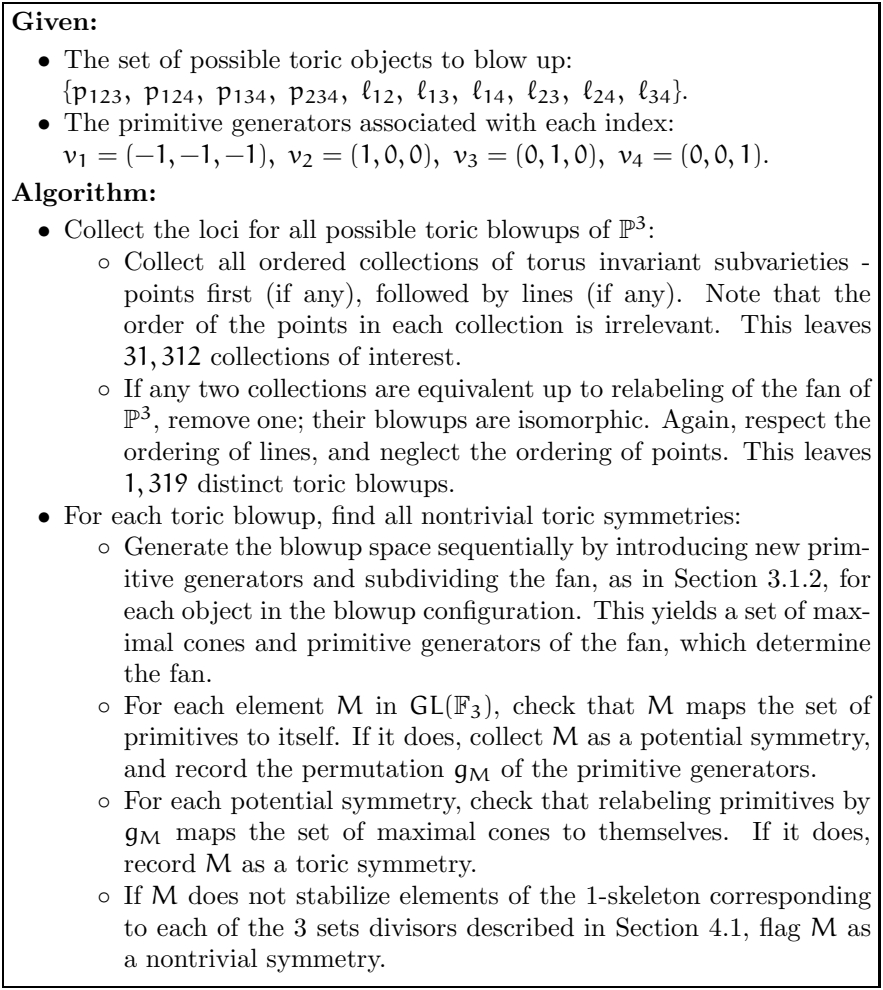

\fbox{\parbox{4.5in}{
\small

{\bf Given:}
\begin{list}{$\bullet$}{\setlength{\leftmargin}{2em}}

	\item The set of possible toric objects to blow up:\\ $\{p_{123},\ p_{124},\ p_{134},\
	p_{234},\ \ell_{12},\ \ell_{13},\ \ell_{14},\ \ell_{23},\ \ell_{24},\ \ell_{34}\}$.

	\item The primitive generators associated with each index:\\
 	$v_1 = (-1,-1,-1),\ v_2 = (1,0,0),\ v_3 = (0,1,0),\ v_4 = (0,0,1)$.

\end{list}

{\bf Algorithm:}
\begin{list}{$\bullet$}{\setlength{\leftmargin}{2em}}

	\item Collect the loci for all possible toric blowups of $\PP^3$: 
	
		\begin{list}{$\circ$}{}
		
			\item Collect all ordered collections of torus invariant subvarieties - points first (if any), followed by
			lines (if any). Note that the order of the points in each collection is irrelevant. This leaves $31,312$ collections
			of interest.

			\item If any two collections are equivalent up to relabeling of the fan of $\PP^3$, remove one;
			their blowups are isomorphic. Again, respect the ordering of lines, and
			neglect the ordering of points. This leaves $1,319$ distinct toric blowups.

		\end{list}

	\item For each toric blowup, find all nontrivial toric symmetries:
		\begin{list}{$\circ$}{}

			\item Generate the blowup space sequentially by introducing new primitive generators
			and subdividing the fan, as in Section~\ref{sec: blowup_configs}, for each object in the
			blowup configuration. This yields a set of maximal cones and primitive
			generators of the fan, which determine the fan.

			\item For each element $M$ in $GL(\mathbb{F}_3)$, check that $M$ maps
			the set of primitives to itself. If it does, collect $M$ as a
			potential symmetry, and record the permutation $g_M$ of the primitive generators.

			\item For each potential symmetry, check that relabeling primitives by $g_M$ maps the
			set of maximal cones to themselves. If it does, record $M$ as a toric symmetry.

			\item If $M$ does not stabilize elements of the 1-skeleton corresponding to each of the 3 sets divisors described in Section~\ref{subsec: nontrivial}, 
			flag $M$ as a nontrivial symmetry.

		\end{list}

\end{list}

}}
\caption{The computational technique by which we exhausted all possible toric symmetries of sequential toric blowups of $\PP^3$}
\label{fig:pseudocode}
\pagebreak
\end{figure}

The result of our search identifies precisely those symmetries of classes A, B, C and D.  This computation also shows that there are no further nontrivial maps that are the pushforwards of toric symmetries, the result of Theorem~\ref{thm: symmetries}. Pseudocode for the computational technique is shown in Figure~\ref{fig:pseudocode}. \qed 

\begin{remark}\label{rem: pf-p2symmetries}
Precisely the same method that characterizes nontrivial symmetries in toric blowups of $\PP^3$ characterizes toric blowups of $\PP^2$, which are blowups of $\PP^2$ at a subset of its three $T$-fixed points. Again, a toric symmetry of $X$, a toric blowup of $\PP^2$, is nontrivial if and only if it exchanges exceptional and non-exceptional classes. Performing the computation yields the results of Theorem~\ref{thm: symmetriesP2}.
\end{remark}

\subsection{Proof of Theorem~\ref{thm: classA}}\label{sec: proof_1}
We will choose a point $p_1$ and lines $\ell_1$ and $\ell_2$ which are fixed by the torus action, and satisfy the intersection conditions of Theorem~\ref{thm: classA}. Since GW and DT invariants are invariant under deformations, proving the result in this case is sufficient. Choose $X$ to be the toric blowup of $\PP^3$ at $p_{123}$, $\ell_{34}$ and $\ell_{24}$, yielding a space of Class A. Consider the action of a lattice isomorphism $\tau$ on $\Sigma_X\subset \ZZ^3$, given by 
\[
\tau = \begin{pmatrix}
0 & 1 & 0\\
1 & 0 & 0\\
1 & 1 & -1
\end{pmatrix}.
\]

\begin{figure}[b]
\def\lineseg{2}
\begin{tikzpicture}[thick,draw=gray,fill=gray,text=gray]
	\tikzstyle{locus}=[draw=red,fill=red,text=red]

	\draw (30:\lineseg) -- (150:\lineseg) node [midway,anchor=south] {$\ell_{14}$};
	\draw[style=locus] (150:\lineseg) -- (270:\lineseg) node [midway,anchor=east] {$\ell_{34}$};
	\draw[style=locus] (270:\lineseg) -- (30:\lineseg) node [midway,anchor=west] {$\ell_{24}$};

	\draw (30:\lineseg) -- (0,0) node [midway,shape=circle,inner sep=1,fill=white,opacity=0.9,text opacity=1] {$\ell_{12}$};
	\draw (150:\lineseg) -- (0,0) node [midway,shape=circle,inner sep=1,fill=white,opacity=0.9,text opacity=1] {$\ell_{13}$};
	\draw (270:\lineseg) -- (0,0) node [midway,shape=circle,inner sep=1,fill=white,opacity=0.9,text opacity=1] {$\ell_{23}$};

	\filldraw[style=locus] (0,0) circle (1pt) node [anchor=south,above=2pt] {$p_{123}$};
	\filldraw (30:\lineseg) circle (1pt) node [anchor=west] {$p_{124}$};
	\filldraw (150:\lineseg) circle (1pt) node [anchor=east] {$p_{134}$};
	\filldraw (270:\lineseg) circle (1pt) node [anchor=north] {$p_{234}$};
\end{tikzpicture}
\caption{\small{The torus-invariant subvarieties of $\PP^3$. The red subvarieties comprise the blowup locus for Class A, blown up in the following order: $p_{123},\ \ell_{34},\ \ell_{24}$.}}
\label{fig: thm2pf-blowup-diagram}
\end{figure}
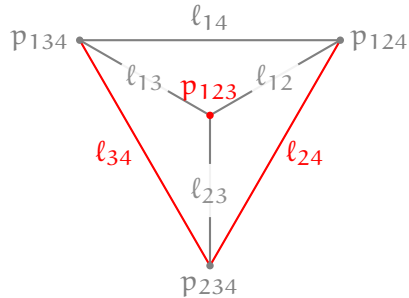

We compute the pushforward $\tau_*$ on $A_*(X)$ using the method described in Section~\ref{subsec: nontrivial} via Poincar\'e duality and obtain the following map on the (basis of) divisor classes $A_2(X) = H_4(X;\ZZ)$. 
\begin{eqnarray*}
\tau_*H &=& 2 H - E_{123} - F_{34}-F_{24}\\
\tau_* E_{123} &=& H - F_{34}-F_{24}\\
\tau_* F_{34} &=& H - E_{123} -F_{24}\\
\tau_* F_{24} &=& H - E_{123} - F_{34}.
\end{eqnarray*}
We compute the intersection product of the Poincar\'e duals of the divisor classes to get the action of $\tau_*$ on $A_1(X)$. The nonzero intersections are given by
\begin{eqnarray*}
H\cdot H &=& h\\
E_{123}\cdot E_{123} &=& -e_{123}\\
F_{34}\cdot F_{34} &=& -f_{34}-s_{34}, \quad s_{34} = h-3f_{34}\\
F_{24}\cdot F_{24} &=& -f_{24}-s_{24}, \quad s_{24} = h-2f_{24}-f_{34}\\
F_{34}\cdot F_{24} &=& f_{24}\\
H\cdot F_{ij} &=& f_{ij},
\end{eqnarray*}
where as above, $h$ is the pullback of the class of a line in $\PP^3$, $e_{123}$ is the class of a line in the exceptional divisor $E_{123}$, $f_{ij}$ is the fiber class in $F_{ij}$ and $s_{ij}$ represents the section class in $F_{ij}$. Note that $F_{ij}$ has trivial fibration and is abstractly isomorphic to $\PP^1\times \PP^1$. The action of $\tau$ on curve classes is given by
\begin{eqnarray*}
\tau_*h &=& 2 h - e_{123} - f_{34}-f_{24}\\
\tau_* e_{123} &=& h - f_{34}-f_{24}\\
\tau_*f_{34} &=& h - e_{123} -f_{24}\\
\tau_*f_{24} &=& h - e_{123} - f_{34}.
\end{eqnarray*}
This completely describes the action of $\tau$ since the classes $h$, $e_{123}$ and $f_{ij}$ form a basis of $H_{2}(X)$. Thus for any $\beta \in H_{2}(X)$, we may write $\beta = dh -a_{123}e_{123}-a_{34}f_{34}-a_{24}f_{24}$ for some integers $d,a_{123},a_{34}$ and $a_{34}$. Therefore by Equations~\ref{eq: gw-symmetry} and~\ref{eq: dt-symmetry}, we have
\[
GW^{X}_{g,\beta}=GW^{X}_{g,\tau_{*}\beta}=GW^{X}_{g,\beta '},
\]
where $\beta '= d'h-a_{123}'e_{123}-a_{24}'e_{24}-a_{34}'e_{34}$ is given in Theorem~\ref{thm: classA}, under the relabeling $a_{123}=a_{1}$, $a_{24}=a_{2}$ and $a_{34}=a_{3}$. \qed

\subsection{Proof of Theorem~\ref{thm: classB}}\label{sec: proof_2}
We will choose points $p_1 = p_{123}$, $p_2 = p_{124}$ and lines $\ell_1 = \ell_{23}$, $\ell_2 = \ell_{34}$ and $\ell_{3} = \ell_{14}$. Now let $X$ be the toric blowup of $\PP^3$ at $p_1$, $p_2$, $\ell_1$, $\ell_2$ and $\ell_3$ yielding a space of Class B. Consider the action of the lattice isomorphism $\tau$ on $\Sigma_X \subset \ZZ^n$, where $\tau$ is given by 
\[
\begin{pmatrix} 0 & 0 & 1 \\ 1 & -1 & 1 \\ 1 & 0 & 0
	\end{pmatrix} .
\]

\begin{figure}[b]
\def\lineseg{2}
\begin{tikzpicture}[thick,draw=gray,fill=gray,text=gray]
	\tikzstyle{locus}=[draw=red,fill=red,text=red]

	\draw[style=locus] (30:\lineseg) -- (150:\lineseg) node [midway,anchor=south] {$\ell_{14}$};
	\draw[style=locus] (150:\lineseg) -- (270:\lineseg) node [midway,anchor=east] {$\ell_{34}$};
	\draw (270:\lineseg) -- (30:\lineseg) node [midway,anchor=west] {$\ell_{24}$};

	\draw (30:\lineseg) -- (0,0) node [midway,shape=circle,inner sep=1,fill=white,opacity=0.9,text opacity=1] {$\ell_{12}$};
	\draw (150:\lineseg) -- (0,0) node [midway,shape=circle,inner sep=1,fill=white,opacity=0.9,text opacity=1] {$\ell_{13}$};
	\draw[style=locus] (270:\lineseg) -- (0,0) node [midway,shape=circle,inner sep=1,fill=white,opacity=0.9,text opacity=1] {$\ell_{23}$};

	\filldraw[style=locus] (0,0) circle (1pt) node [anchor=south,above=2pt] {$p_{123}$};
	\filldraw[style=locus] (30:\lineseg) circle (1pt) node [anchor=west] {$p_{124}$};
	\filldraw (150:\lineseg) circle (1pt) node [anchor=east] {$p_{134}$};
	\filldraw (270:\lineseg) circle (1pt) node [anchor=north] {$p_{234}$};
\end{tikzpicture}
\caption{\small{The torus-invariant subvarieties of $\PP^3$. The red subvarieties comprise the blowup locus for Class B, blown up in the following order: $p_{123},\ p_{124},\ \ell_{23},\ \ell_{34},\ \ell_{14}$.}}
\label{fig: thm3pf-blowup-diagram}
\end{figure}
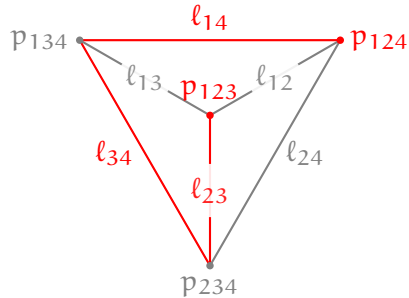

Pushing forward $\tau$ to $A_2(X)$ as above, we compute 
\begin{eqnarray*}
\tau_* H &=& 2 H - E_{123} - E_{124} - F_{23}-F_{34}\\
\tau_* E_{123} &=& F_{14}\\
\tau_* E_{124} &=& H - E_{123} -  F_{23}-F_{34}\\
\tau_* F_{23} &=& H-E_{124}-F_{34}-F_{14}\\
\tau_* F_{34} &=& H -E_{123}-E_{124}-F_{23}\\
\tau_* F_{14} &=& E_{123}.
\end{eqnarray*}
The nonzero intersection pairings are as follows.
\begin{eqnarray*}
H\cdot H &=& h\\
E_{ijk}\cdot E_{ijk} &=& -e_{ijk}\\
F_{ij}\cdot F_{ij} &=& -f_{ij}-s_{ij}, \quad s_{23} = h-e_{123}-f_{23}\\
\phantom{F_{ij}\cdot F_{ij}} &\phantom{=}& \phantom{-f_{ij}-s_{ij},}\,\quad s_{34} = h-f_{23}-2f_{34}\\
\phantom{F_{ij}\cdot F_{ij}} &\phantom{=}& \phantom{-f_{ij}-s_{ij},}\,\quad s_{14} = h-e_{124}-f_{34}\\
F_{34}\cdot F_{14} &=& f_{14}\\
F_{23}\cdot F_{34} &=& f_{34}\\
E_{ijk}\cdot F_{ij} &=& f_{ij}\\
H\cdot F_{ij} &=& f_{ij}.
\end{eqnarray*}
The action of $\tau$ on curve classes is given by
\begin{eqnarray*}
\tau_* h &=& 2 h - e_{124} - f_{23}-f_{34}\\
\tau_* e_{123} &=& h-e_{124}-f_{34}+ f_{14}\\
\tau_* e_{124} &=& h - f_{23}-f_{34}\\
\tau_* f_{23} &=& h-e_{124}-f_{34}\\
\tau_* f_{34} &=& h - e_{124}-f_{23}\\
\tau_* f_{14} &=& e_{123} - f_{23}.
\end{eqnarray*}
Finally, applying $\tau_*$ to an arbitrary curve class $\beta\in H_2(X)$, we obtain Theorem~\ref{thm: classB}. \qed 

\subsection{Proof of Theorem~\ref{thm: classC}}\label{sec: proof_3}
We proceed as in the proofs of Theorems~\ref{thm: classA} and~\ref{thm: classB}. Choose points $p_1 = p_{124}$, $p_2 = p_{123}$ and lines $\ell_1 = \ell_{34}$, $\ell_2 = \ell_{23}$ and $\ell_{3} = \ell_{14}$. Let $X$ be the toric blowup of $\PP^3$ at $p_1$, $p_2$, $\ell_1$, $\ell_2$ and $\ell_3$, yielding a space of Class C. Consider the action of the lattice isomorphism $\tau$ on $\Sigma_X \subset \ZZ^n$, where $\tau$ is given by 
\[
\tau = \begin{pmatrix} 0 & 0 & -1 \\ 1 & 0 & -1 \\ 1 & -1 & 0
	\end{pmatrix}
\]

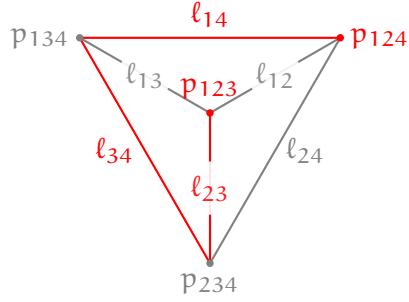
\begin{figure}
\def\lineseg{2}
\begin{tikzpicture}[thick,draw=gray,fill=gray,text=gray]
	\tikzstyle{locus}=[draw=red,fill=red,text=red]

	\draw[style=locus] (30:\lineseg) -- (150:\lineseg) node [midway,anchor=south] {$\ell_{14}$};
	\draw[style=locus] (150:\lineseg) -- (270:\lineseg) node [midway,anchor=east] {$\ell_{34}$};
	\draw (270:\lineseg) -- (30:\lineseg) node [midway,anchor=west] {$\ell_{24}$};

	\draw (30:\lineseg) -- (0,0) node [midway,shape=circle,inner sep=1,fill=white,opacity=0.9,text opacity=1] {$\ell_{12}$};
	\draw (150:\lineseg) -- (0,0) node [midway,shape=circle,inner sep=1,fill=white,opacity=0.9,text opacity=1] {$\ell_{13}$};
	\draw[style=locus] (270:\lineseg) -- (0,0) node [midway,shape=circle,inner sep=1,fill=white,opacity=0.9,text opacity=1] {$\ell_{23}$};

	\filldraw[style=locus] (0,0) circle (1pt) node [anchor=south,above=2pt] {$p_{123}$};
	\filldraw[style=locus] (30:\lineseg) circle (1pt) node [anchor=west] {$p_{124}$};
	\filldraw (150:\lineseg) circle (1pt) node [anchor=east] {$p_{134}$};
	\filldraw (270:\lineseg) circle (1pt) node [anchor=north] {$p_{234}$};
\end{tikzpicture}
\caption{\small{The torus-invariant subvarieties of $\PP^3$. The red subvarieties comprise the blowup locus for Class C, blown up in the following order: $p_{124},\ p_{123},\ \ell_{34},\ \ell_{23},\ \ell_{14}$.}}
\label{fig: thm4pf-blowup-diagram}
\end{figure}

Pushing forward $\tau$ to $H_4(X)$, we see that 
\begin{eqnarray*}
\tau_* H &=& 2 H - E_{124} - E_{123} - F_{34}-F_{14}\\
\tau_* E_{124} &=& H - E_{124} -  F_{34}-F_{14}\\
\tau_* E_{123} &=& F_{23}\\
\tau_* F_{34} &=& H-E_{124}-E_{123}-F_{14}\\
\tau_* F_{23} &=& H - E_{123}-F_{34}-F_{23}\\
\tau_* F_{14} &=& E_{124}.
\end{eqnarray*}
Now, computing the intersection of the divisor classes, we find the nonzero intersections to be:
\begin{eqnarray*}
H\cdot H &=& h\\
E_{ijk}\cdot E_{ijk} &=& -e_{ijk}\\
F_{ij}\cdot F_{ij} &=& -f_{ij}-s_{ij}, \quad s_{34} = h-3f_{34}\\
\phantom{F_{ij}\cdot F_{ij}} &\phantom{=}& \phantom{-f_{ij}-s_{ij},}\,\quad s_{23} = h-e_{123}-f_{34}\\
\phantom{F_{ij}\cdot F_{ij}} &\phantom{=}& \phantom{-f_{ij}-s_{ij},}\,\quad s_{14} = h-e_{124}-f_{34}\\
F_{34}\cdot F_{23} &=& f_{23}\\
F_{34}\cdot F_{14} &=& f_{14}\\
E_{ijk}\cdot F_{ij} &=& f_{ij}\\
H\cdot F_{ij} &=& f_{ij}.
\end{eqnarray*}
Using the above intersections, again, we compute the map $\tau_*$ on the basis of $A_1(X)$:
\begin{eqnarray*}
\tau_* h &=& 2 h - e_{123} - f_{34}-f_{14}\\
\tau_* e_{124} &=& h - f_{34}-f_{14}\\
\tau_* e_{123} &=& h-e_{123}-f_{34}+ f_{23}\\
\tau_* f_{34} &=& h-e_{123}-f_{14}\\
\tau_* f_{23} &=& h - e_{123}-f_{34}\\
\tau_* f_{14} &=& e_{124} - f_{14}.
\end{eqnarray*}
Now similarly, choose the matrix $\sigma$ below, and push-forward $\sigma$ to $A_*(X)$. 
\[
\sigma = \begin{pmatrix} 1 & -1 & 0 \\ 0 & -1 & 0 \\ 0 & 0 & -1
	\end{pmatrix}.
\]
This yields the following map on the basis of $A_1(X)$,
\begin{eqnarray*}
\sigma_* h &=& 2 h - e_{123} - f_{34}-f_{14}\\
\sigma_* e_{124} &=& h - e_{123}-f_{34}+f_{23}\\
\sigma_* e_{123} &=& h-f_{34}-f_{14}\\
\sigma_* f_{34} &=& h-e_{123}-f_{14}\\
\sigma_* f_{23} &=& e_{124}-f_{14}\\
\sigma_* f_{14} &=& h-e_{123}-f_{34}.
\end{eqnarray*}
Applying $\sigma_*$ and $\tau_*$ to an arbitrary element $\beta\in A_1(X)$, we obtain Theorem~\ref{thm: classC}. Additionally, note that the set $\{\tau,\tau^2,\sigma,\sigma\tau\}$ are precisely the four nontrivial symmetries of $X$ found in Theorem~\ref{thm: symmetries}. \qed

\subsection{Higher Virtual Dimension}\label{subsec: nonzero_vd} The results stated in Section~\ref{sec: intro} appear without insertions of cohomology classes for the sake of brevity. However, the results hold in higher virtual dimension as well. 

Let $X$ be a variety of Class A, B, C or D, and let $\tau : X\rightarrow X$ be the automorphism of $X$ induced via toric symmetry. Now let $\beta \in H_{2}(X)$ and $\gamma'_1,\ldots,\gamma'_r \in H^*(X)$. Then, by Equations~\ref{eq: gw-symmetry} and~\ref{eq: dt-symmetry} and the proofs of Theorems 2 through 5, we have
\begin{equation}\label{eq: GW-insertions}
GW^{X}_{n,\tau_{*}\beta}(\gamma_{1},\dotsc ,\gamma_{r}) =
GW^{X}_{n,\beta}(\tau^{*} \gamma_{1},\dotsc ,\tau^{*}\gamma_{r}),
\end{equation}
and
\begin{equation}\label{eq: DT-insertions}
DT^{X}_{n,\tau_{*}\beta}(\gamma_{1},\dotsc ,\gamma_{r}) =
DT^{X}_{n,\beta}(\tau^{*} \gamma_{1},\dotsc ,\tau^{*}\gamma_{r}),
\end{equation}
where the action of $\tau$ on curve classes and cohomology is given above in Sections~\ref{sec: proof_1}--\ref{sec: proof_3}.

\bibliographystyle{plain}
\bibliography{p3symmetries}

\end{document}